\documentclass[10pt,twocolumn,twoside]{article}

\usepackage{palatino}  
\usepackage[T1]{fontenc}

\usepackage[left=3.5cm,top=2.54cm,right=2.5cm,bottom=2.8cm,nohead]{geometry} 
\usepackage{amsmath}
\usepackage[font=sf]{caption}
\usepackage{color}
\usepackage{fancyhdr}
\usepackage{graphicx}
\usepackage{lettrine}
\usepackage[super,sort&compress]{natbib}
\usepackage{pifont}
\usepackage{sectsty}
\usepackage{setspace}
\usepackage{sidecap}  
\usepackage{subfig}
\usepackage{wrapfig}

\definecolor{sectcol}{rgb}{0.0,0.24,0.43} 
\definecolor{dropped}{rgb}{0.55,0.06,0.11}
\definecolor{update}{rgb}{0,0,0} 

\sectionfont{\large\sffamily\color{sectcol}\vspace{-2mm}}
\subsectionfont{\normalsize\sffamily\bfseries\vspace{-2mm}}

\newlength{\up}
\usepackage[
    pdftex,
    pdftitle={},
    pdfauthor={Rio Yokota, L.~A. Barba},
    pdfpagemode={UseOutlines},
    bookmarks, bookmarksopen,bookmarksnumbered={True},
    colorlinks, linkcolor={black},citecolor={black},urlcolor={black}
    ]{hyperref}
    
\begin{document}
\pagestyle{fancy}

\twocolumn[ 
{\Huge{\color{sectcol}Hierarchical $N$-body simulations  with auto-tuning for heterogeneous systems}}
\vspace{0.8cm}

{ \sf 
\onehalfspacing
{\large \textit{Algorithms designed to efficiently solve this classical problem of physics  fit very well on GPU hardware, and exhibit excellent scalability on many GPUs. Their computational intensity makes them a promising approach for many other applications amenable to an N-body formulation. Adding features such as auto-tuning makes multipole-type algorithms ideal for heterogeneous computing environments. }}
\vspace{1cm}

\textbf {Rio Yokota, Lorena A. Barba}\\
  Mechanical Engineering Dept., Boston University, Boston MA 02215

\vspace{1cm}
}
] 

\lettrine{\textcolor{dropped}{T}}{} he classic $N$-body problem of mechanics solves for the motion of $N$ bodies interacting via the force of gravitation. Beyond gravitational masses, a variety of physical systems can be modeled by the interaction of $N$ particles, e.g., atoms or ions under electrostatics and van der Waals forces lead to molecular dynamics. Also, the integral formulation of problems modeled by elliptic partial differential equations leads to numerical integration having the same form, computationally, as an $N$-body interaction. \textcolor{update}{Astrophysics and molecular dynamics problems rely on the Laplace potentials, while Helmholtz potentials have applications in acoustics and electromagnetics, and Stokes potentials can be found in elasticity and geophysics.} Adding to this diversity of applications, radiosity algorithms for global illumination problems in computer graphics also benefit from $N$-body methods.

In the absence of a closed-form solution beyond 3 bodies, $N$-body problems require a numerical approach. The direct simulation of the \emph{all-pairs} interaction results in a computational complexity of order $N^2$, which becomes too expensive to compute for large $N$. Nevertheless, the simplicity of direct integration admits ease of use of hardware accelerators leading to a prominent branch of research in this area. Beyond this approach, many notable algorithmic inventions bear on fast computation of $N$-body interactions. Among them, multipole-based methods are gaining traction as ideally placed for the heterogeneous, many-core hardware environment emerging beyond the petascale computing era. As we show below, fast algorithms such as treecodes and the fast multipole method are characterized by a high computational intensity, as measured using the roofline model. This reflects on their excellent performance on GPU hardware.

\section*{History of direct \textit{N}-body simulation \& special-purpose machines}
\vspace{\up}

Numerical simulation of many-body dynamics on digital computers began in the early 60's in two distinct fields of physics: astrophysics and molecular dynamics. These first simulations were able to compute only in the order of a $100$ particles.  Since then, the $N$-body community has done so much more than simply rely on half a century's worth of Moore's law-governed performance improvements. Considerable effort has been dedicated to algorithmic innovations, special-purpose hardware, and performance optimizations. 

The $N$-body problem of astrophysics was such a strong motivator in computational science, that it drove the creation of a special-purpose supercomputer consisting of dedicated pipelines for $N$-body simulations. GRAPE-1 was the first of its kind and achieved a performance comparable to the CRAY-XMP/1 at 1/10,000 the cost. The $4^{\text{th}}$ generation GRAPE-4 was the first computer to reach 1 teraflop/s, but without being able to perform the LINPACK benchmark, it could not be so recognized in the Top500 list. The GRAPE machines were dedicated to gravitational $N$-body simulations, but were later extended to perform molecular dynamics computations and given the name MDGRAPE. However, as the cost of fabricating a micro-processor became exceedingly high, it became difficult for a single research group to produce its own processor. This situation was exacerbated by the arrival of GPUs, which by providing the same capability at much lower cost have emerged as a disruptive technology for $N$-body simulations.

\section*{Gordon Bell prize record}
\vspace{\up}

$N$-body simulations have persistently been at the forefront of high-performance computing in terms of the flop/s that they deliver. Groundbreaking $N$-body simulations have won the coveted Gordon Bell prize of computing fourteen times from 1992 to 2010. 
Some may argue that it is not the achieved flop/s that counts, but the science that they deliver. This is a valid assertion, in the same sense that reaching exaflop/s in itself should not be the purpose of high-performance computing. The maximum flop/s for $N$-body simulations is obtained with the pleasingly parallel $\mathcal{O}(N^2)$ all-pairs summation, whereas the maximum amount of science will often be delivered by a fast algorithm. Many of the award-winning $N$-body simulations have used hierarchical $N$-body algorithms, discussed next, and not the all-pairs summation.

\section*{Hierarchical \textit{N}-body algorithms}

\subsection*{Introduction to treecodes \& FMM}
\vspace{\up}

The amount of computation required by a direct $N$-body simulation is $\mathcal{O}(N^2)$, which quickly becomes prohibitive for large $N$, even with special-purpose machines. The treecode algorithm\cite{BarnesHut1986}  reduces the complexity to $\mathcal{O}(N\log N)$ by clustering the remote particles into progressively larger groups and using multipole expansions to approximate their influence on each target particle. The fast multipole method\cite{GreengardRokhlin1987}, FMM, can achieve $\mathcal{O}(N)$ by clustering not only the remote particles, but also the nearby particles using local expansions. \textcolor{update}{We give a technical but brief overview of the two algorithms in our chapter in the latest \textit{GPU Gems} volume\cite{YokotaBarba2010}.}
Both treecodes and FMM use tree data structures to cluster particles into a hierarchy of cells. Historically, however, the two methods have followed separate paths of evolution, and have adopted  different practices to meet the demands of their communities of users and their applications.

Treecodes have been popular predominantly in the astrophysics community, where the particle distribution is always highly non-uniform. Thus, the treecode evolved as an inherently adaptive algorithm; on the other hand, the FMM community viewed adaptivity as an additional feature. FMMs have not necessarily focused on a specific application, and practitioners were often aiming at higher accuracy than in the case of treecodes. \textcolor{update}{Typically, the series expansion truncation level can be $p=10-15$, with for example $p=10$ giving an accuracy of 4 significant digits in the potential for the Laplace kernel.} As greater accuracy is expressed by a higher truncation level for the series, a variety of fast translation methods have been developed for FMM (based on spherical and plane wave expansions) that can achieve higher accuracy at reduced cost (e.g., $p^{4}$ rather than $p^{6}$). Most treecodes, in contrast, use simple Taylor series of low order in Cartesian coordinates; \textcolor{update}{typically, $p=3$ is used in practice}. Treecodes also use the ratio between the size of cells $b$ and the distance $l$ between them to construct the interaction list. This is known as the multipole acceptance criterion (MAC), \textcolor{update}{$\theta=b/l$, and it is used to determine if a cell should be evaluated or subdivided further; a smaller value will increase the accuracy of the approximation.} In contrast, FMMs use parent, child, and neighbor relationships to construct the interaction list. These differences between treecodes and FMM are mostly due to historical reasons, rather than mathematical or algorithmic ones. Hence, cross-fertilization of these two fields is surely possible, and may produce an algorithm that takes advantage of the best features of both methods.

\subsection*{Hybrid treecode/FMM algorithm}
\vspace{\up}

One main difference between treecodes and FMM is the fact that treecodes calculate cell-particle interactions, while FMM calculates cell-cell interactions for the far field. Warren \& Salmon \cite{WarrenSalmon1995} suggested a technique to calculate cell-cell interactions in treecodes, but very little emphasis was placed on this technique in their paper. That work was elucidated and refined by Dehnen \cite{Dehnen2002}, who introduced other techniques such as generic tree traversals, mutual cell-cell interactions, and error-controlled multipole acceptance criteria.

From the FMM side, Cheng \textit{et al.}\cite{ChengETal1999} discussed a mechanism to select between cell-particle and cell-cell interactions that, according to the authors, \textcolor{update}{should always be faster} than a pure treecode or pure FMM. They also developed a more adaptive cell-cell interaction stencil that considers the interaction of cells at different levels in the tree; this is similar to what could be obtained from using the MAC in FMMs, which would naturally allow interaction between different levels by relying on the ratio between  size and distance of cells.

The two hybrid ideas mentioned above---the $\mathcal{O}(N)$ treecode by Dehnen, and the adaptive FMM by Cheng \textit{et al.}---were compared on the same machine by Dehnen \cite{Dehnen2002}. The relative performance depends on the required accuracy, where the $\mathcal{O}(N)$ treecode outperforms the adaptive FMM by an order of magnitude for 4 significant digits of accuracy, while the FMM performs better if the required accuracy is over 6 digits. Dehnen attributes the inefficiency of his code at higher accuracies to the fact that the order of expansions is kept constant while accuracy control is effected using the MAC. This suggests that adding the capability to handle variable order of expansions in Dehnen's framework could produce a very fast treecode-FMM hybrid method.

\medskip

We have recently developed a hybrid treecode-FMM that has similar structure to Dehnen's method but has control over both the order of expansion and MAC. Our kernels are based on spherical harmonic expansions, but have the capability to switch to Cartesian expansions if the required accuracy is lower than a certain threshold; this is the key to achieving high performance for low-accuracy calculations. In addition, we have incorporated a key feature for achieving high performance in today's hardware:  the capability to auto-tune the kernels on heterogenous architectures; we will explain this feature in detail below.

\section*{Use of GPUs for \textit{N}-body simulation}

\subsection*{Early application of GPUs}
\vspace{\up}

When CUDA 1.0 was released in 2007 and graphics cards became programmable in C, there were very few scientific applications that could take advantage of this new programming paradigm. $N$-body simulations were one of the first applications to extract the full compute capability of GPUs. Hamada \& Iitaka's initial effort\cite{HamadaIitaka2007}  parallelized the \textit{source} particles among thread blocks, which required a large reduction to be performed at the end. Soon after that, Nyland~\textit{et al.}\ published their chapter\cite{NylandETal2007} in the \textsl{GPU Gems 3} book, using the opposite approach: the \textit{target} particles were parallelized among thread blocks and no reduction was necessary. Relying on the same technique, Belleman~\textit{et al.}\ released  in 2008 their code named \textit{Kirin} (after a Japanese beer) \cite{BellemanETal2008}. In 2009, Gaburov~\textit{et al.}\ emulated the GRAPE-6 machine in their GPU code \textit{Sapporo} (named after another Japanese beer) \cite{GaburovETal2009}, which they were able to do thanks to the similarities between the GRAPE and GPU architectures. The fact that special-purpose GRAPEs were so similar to GPUs may have given the $N$-body community an advantage, since many techniques that were developed a decade ago to tune the codes for GRAPE could be used directly on GPUs.

\subsection*{Advantage of \textit{N}-body algorithms on GPUs}
\vspace{\up}

Perhaps it is not a coincidence that current GPUs turn out to have similarities to the GRAPE hardware. Computation did not suddenly become cheap, and communication did not suddenly become comparatively more expensive. The trend has always been there, and data-parallel architectures like GRAPE had been performing much better than serial processors all along. It was only a matter of time until mass-produced data-parallel processors appeared and took over certain application areas. In retrospect, the late 90's and early 2000's were peculiar times in the history of high-performance computing, when even the most data-parallel algorithms were being computed on serial processors. Note that we refer to fine-grained parallelism, and not the coarse-gained parallelism that \textit{was} properly handled in parallel during this era with MPI.

\begin{figure}
\centering
\includegraphics[width=0.49\textwidth]{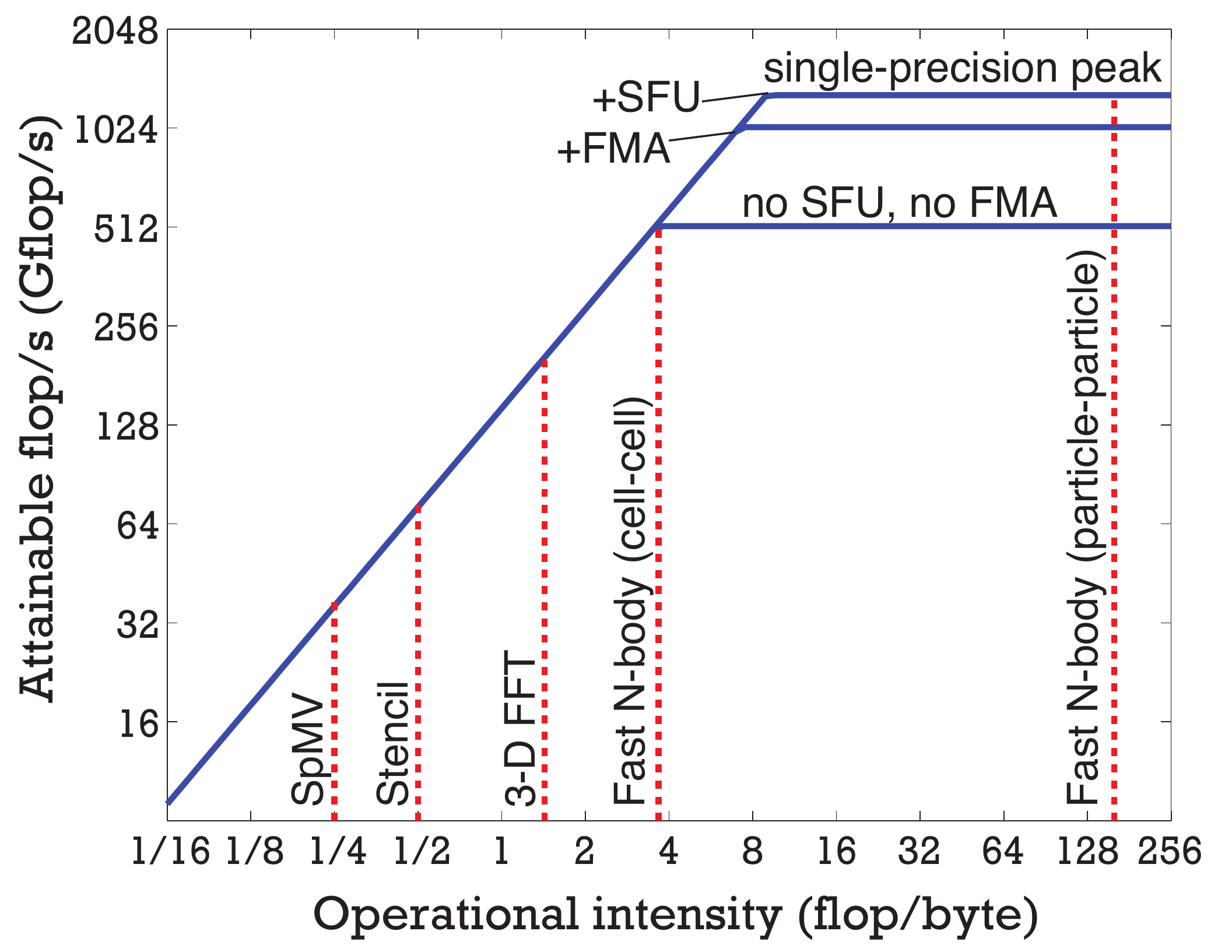}
\caption{Roofline model of FMM kernels on an NVIDIA C2050 GPU. The `SFU' label is used to indicate the use of special function units and `FMA' indicates the use of fused multiply-add instructions. The order of multipole expansions was set to $p=15$.}
\label{fig:roofline}
\end{figure}

We quantify the advantage of $N$-body algorithms on GPUs via the roofline model \cite{WilliamsETal2008}. This model offers a useful metric for predicting the performance of algorithms on multicore architectures, with the number of floating-point operations per byte of data transferred used to determine whether the algorithm will be limited by the floating-point performance of the processor, or by the memory bandwidth.
The roofline model for any algorithm is contingent on the hardware used.  We show the roofline on an NVIDIA Tesla C2050 GPU in Figure \ref{fig:roofline}, including both the particle-particle and cell-cell interactions of the FMM. The C2050 has a memory bandwidth of 144 GB/s and a single-precision peak performance of 1288 Gflop/s when special function units (SFU) and fused multiply-add (FMA) operations are fully utilized. Without the use of SFUs, the peak performance decreases to 1030.4 Gflop/s, and further down to 515.2 Gflop/s if the FMA is not applicable. Fortunately, $N$-body kernels have many adjacent multiply-add operations that can be fused.

As seen on Figure~\ref{fig:roofline}, the particle-particle interaction is pleasingly parallel, with the operational intensity reaching 160 and well under the flat part of the roofline. The operational intensity of the cell-cell interaction is not as high, but it is still much higher than most algorithms. Figure~\ref{fig:roofline} includes a sparse matrix-vector multiplication (labelled `SpMV'), a multigrid method with a seven-point stencil (`Stencil'), and a 3D fast Fourier transform (`3D FFT') under the same roofline\cite{WilliamsETal2008}. In summary, the roofline model distinctly quantifies the high operational intensity of fast $N$-body algorithms, and reveals their unmistakable advantage on many-core architectures.

\subsection*{Domain decomposition in fast $N$-body methods}
\vspace{\up}

Multi-GPU implementations are indispensable for solving large-scale $N$-body problems, and traditional MPI-based parallelization must be combined with the GPU kernels in order to achieve this. The key to successfully parallelizing fast $N$-body algorithms on distributed memory architectures is the partitioning and communication of the tree structure among individual processes. Salmon \& Warren\cite{Salmon1990,WarrenSalmon1993} made a significant contribution to this area by introducing techniques such as orthogonal recursive bisection (ORB), the local essential tree (LET), and $N$-D hypercube communication. Dubinski \cite{Dubinski1996} summarizes their efforts in a concise and clear manner. A recent improvement in this area is the balancing of linear octrees by Sundar \textit{et al.}\cite{SundarETal2008}, a technique to repartition the domain so that the per-processor partitions are  better aligned with cell boundaries at coarser levels of the octree. This was a key technology behind the 2010 Gordon Bell prize-winning paper \cite{RahimianETal2010}. 

\section*{Many-GPU calculations with FMM}

\subsection*{Biomolecular electrostatics}
\vspace{\up}

We demonstrated the application of the FMM algorithm in a multi-GPU system for biological applications in Yokota \emph{et al.}\cite{YokotaETal2011a} The FMM was used to accelerate a boundary element method solution of the continuum electrostatic model, a popular model for calculating electrostatic interactions between biological molecules in solution. Through guest access to the DEGIMA cluster at Nagasaki University (which holds the \#3 spot in the Green500 list), we were able to test the parallel FMM  on hundreds of GPUs. The largest calculation solved a system of over a billion boundary unknowns for more than 20 million atoms, requiring one minute of run time on 512 GPUs. This work demonstrates that the FMM on GPU could enable routine calculations that were unfeasible before, for example, for analyses of protein-protein interactions in vital biological processes.

\subsection*{Fluid turbulence simulations}
\vspace{\up}

Recently, we applied our periodic FMM algorithm to the simulation of homogeneous turbulent flow in a cube, demonstrating scalable computations with many GPUs.  These calculations were carried out in the TSUBAME 2.0 system, thanks to guest access provided by the Grand Challenge Program of TSUBAME. A total of 2048 NVIDIA M2050 GPUs were used, corresponding to half of the complete system, to achieve a sustained performance of 0.5 petaflop/s. The peak performance of the complete  system is 2.4 petaflop/s, thus we achieved excellent usage of the computational resource in an application.

The preferred method for the simulation of homogeneous isotropic turbulence in a periodic cube has always been the pseudo-spectral method. We compared an FMM-based vortex method with an FFT-based pseudo-spectral method for  turbulence at $Re_{\lambda}\approx 500$ using $2048^3$ grid points, confirming that relevant statistics quantitatively match. The parallel scalability of the FMM algorithm is excellent, obtaining 72\% parallel efficiency in a weak scaling test up to 2048 GPUs.
 With this recent work, we show that the scalability of this algorithm starts to become an advantage over FFT-based methods beyond 2000 parallel MPI processes.
 
 \begin{figure}
\centering
\includegraphics[width=0.40\textwidth]{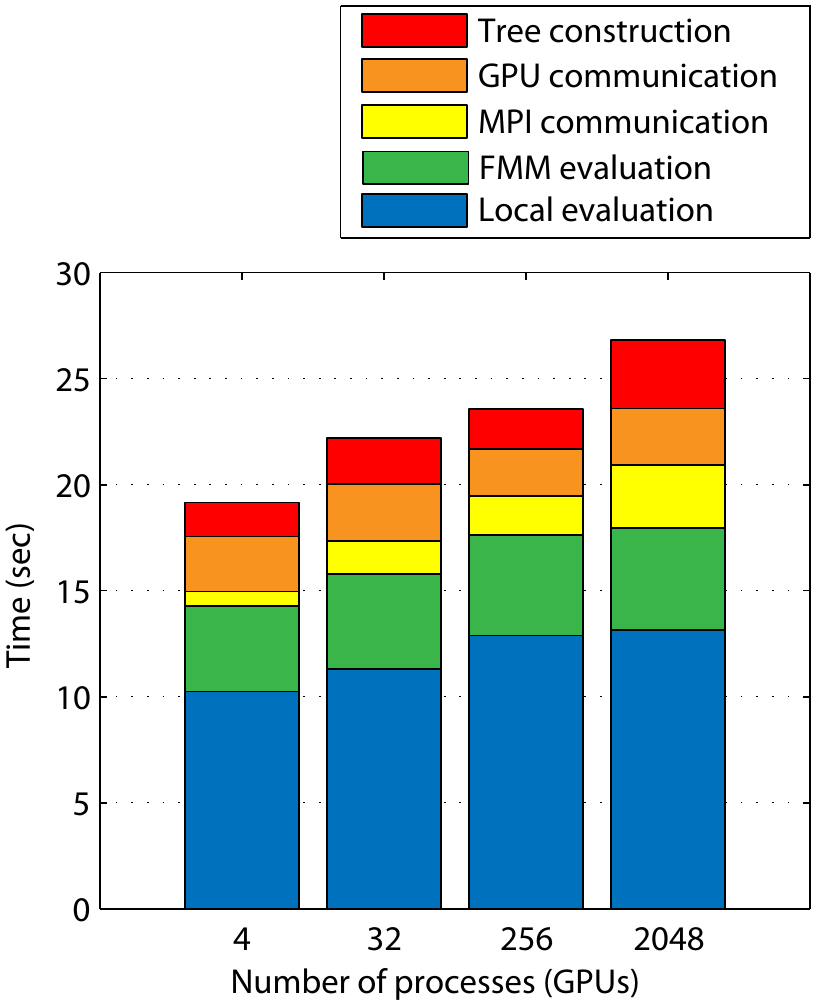}
\caption{Weak scaling test with our FMM code on many GPUs. The parallel efficiency with 2048 GPUs is 72\%.}
\label{fig:weak}
\end{figure}

The results of a weak scaling test with 4 million particles per process is shown in Figure \ref{fig:weak}. The label `Local evaluation' corresponds to the particle-particle kernel, while the `FMM evaluation' label corresponds to the sum of all the other kernel evaluations. The MPI communication is overlapped with the kernel evaluations, so in the bar plot we show the actual time for communications, and the excess time required for the FMM evaluation obtained by subtracting the MPI communications time to the total evaluation time. In this way, the total height of the bar correctly represents the total wall-clock time of the full overlapped calculation.

\section*{A new hybrid treecode-FMM with auto-tuning}

\subsection*{The purpose of auto-tuning}
\vspace{\up}

Like most algorithms, hierarchichal $N$-body methods permit a wide variety of mathematical formulations and computational implementations, some of which are better suited for a particular architecture and not for others. Some of the available choices are: use of Cartesian vs.\ spherical expansions, rotation-based vs.\ plane wave-based translations, cell-cell vs.\ cell-particle interactions for the far field, and choice of order of expansion vs.\ MAC-based error optimization. For each of these choices, the option that will be most efficient depends on (i) the required accuracy, (ii) the hardware and, unfortunately, (iii) the implementation/tuning of kernels. The hardware dependence is particularly problematic when heterogeneous architectures come into the equation. The reason is that the various algorithmic kernels (e.g., cell-cell and cell-particle interaction) achieve different levels of performance measured in flop/s on different hardware. Thus, a well-balanced FMM calculation on one type of hardware might be unbalanced with the same parameters when moving to a different hardware. A simple solution to this problem would be to time all kernels on each hardware, and use this information to select the optimal combination during runtime. For example, if the GPU can perform cell-particle interactions faster than cell-cell interactions for a certain number of particles per cell, our hybrid treecode-FMM will shift more towards treecodes automatically.

\subsection*{Kernel pre-calculation}
\vspace{\up}

The key for auto-tuning in our hybrid fast $N$-body method is the pre-calculation of the kernels. All kernels are evaluated using artificial coordinates, mass/charges, multipole coefficients, and their execution time is measured. This information is then used to select the optimum kernel during the dual tree traversal (described below), choosing between cell-cell, cell-particle, and particle-particle interactions. Certain kernels can achieve higher flop/s than others, which the initial timings will reflect, so that the selection of kernels will be optimized for the architecture automatically. Therefore, any manual parameter tuning associated with hybridizing treecodes and FMMs is rendered unnecessary.

\subsection*{Dual tree traversal}
\vspace{\up}

\begin{figure*}
\centering
\includegraphics[width=0.88\textwidth]{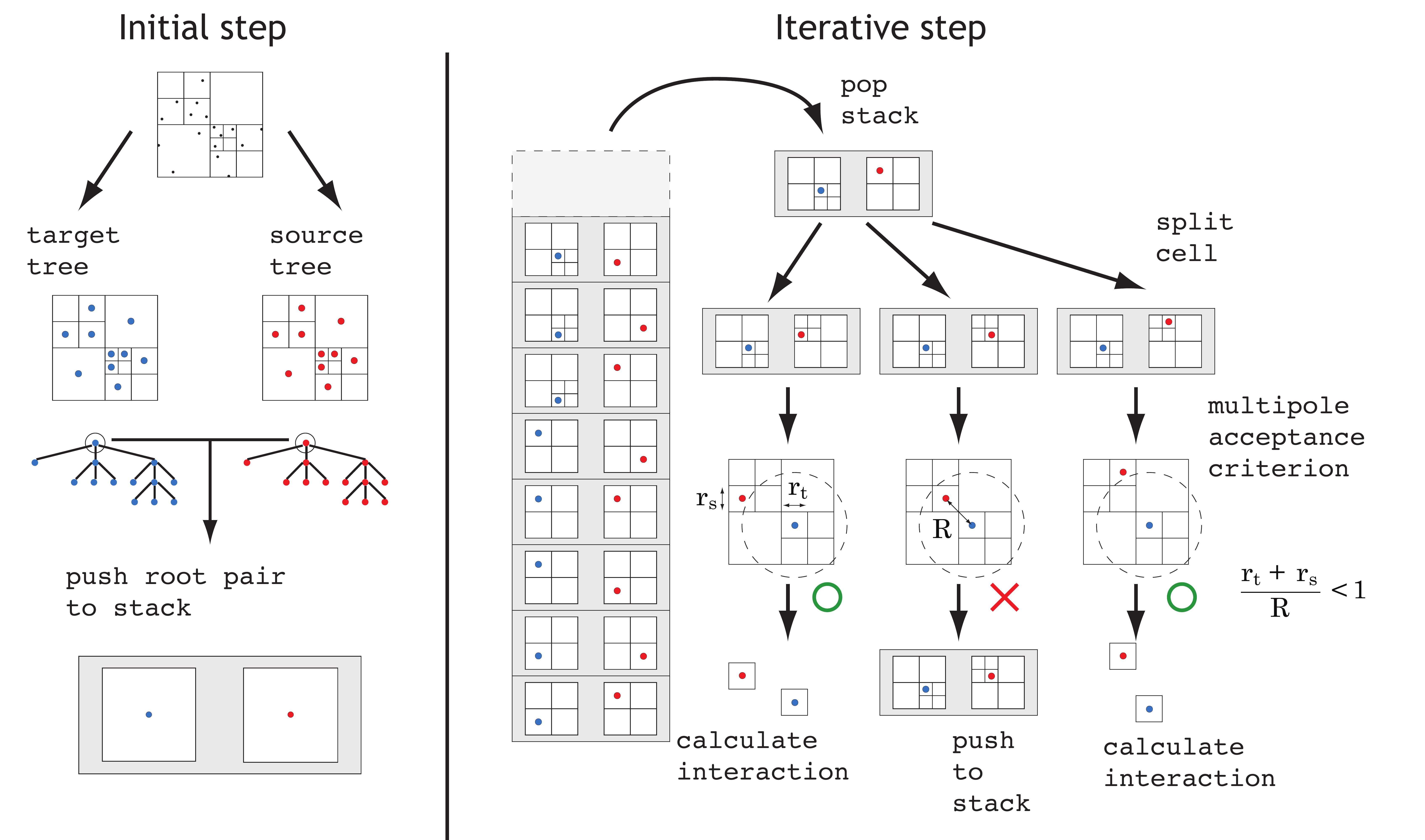}
\caption{Dual tree traversal: illustration of the stack-based procedure.}
\label{fig:treewalk}
\end{figure*}

Auto-tuning of a treecode-FMM hybrid method by dynamically selecting kernels during runtime requires a generic and flexible $\mathcal{O}(N)$ algorithm for traversing the tree. We now describe such a generic tree traversal algorithm, illustrated in Figure \ref{fig:treewalk}. This algorithm can be viewed as a dual tree traversal for the target tree and source tree, which results in $\mathcal{O}(N)$ complexity. In most cases the targets and sources are identical, but the present method can also handle cases where they are not.

The dual tree traversal uses a typical ``last in, first out'' stack data structure that holds pairs of cells: one target cell and one source cell.
As shown on the left panel of Figure \ref{fig:treewalk}, once the tree is constructed, the pair of root cells is pushed into an empty stack. 
After this initial step, the following iterative procedure is applied---as illustrated on the right side of Figure \ref{fig:treewalk}. 
First, a pair of cells is popped from the stack, and the larger of the two cells is subdivided. Always splitting the larger of the two cells guarantees that the pairs in the stack consist of cells of somewhat similar size, which is a necessary condition for achieving $\mathcal{O}(N)$. 
Let us assume, without loss of generality, that the source cell was subdivided (as shown in the Figure). Next, its offspring are matched with the target cell to form new pairs. If a matching set is composed of leaf cells (at the terminal level of the tree), a direct summation is performed between all particles in the cells. If not, the multipole acceptance criterion is used to examine each of the newly created pairs of cells. 
If the cells in a new pair are far/small enough, the interaction between the two cells is immediately calculated. The type of interaction---cell-cell vs.\ cell-particle, rotation-based vs.\ plane wave-based translation---is then chosen to optimize the performance. If the cells are too close/large, then the pair of cells is pushed to the top of the stack. This procedure then starts again and is repeated until the stack is empty.

The procedure described is a simple but highly adaptive and flexible way of performing an $\mathcal{O}(N)$ tree traversal. Interestingly, traditional FMMs do not rely on such algorithms; they construct instead a rigid interaction list for every target cell using parent, child, and neighbor relationships. This in turn requires the kinship in the tree to be directly associated to the geometrical proximity of the cells, i.e., all cells must be perfect cubes. In contrast, the dual tree traversal can be applied to adaptive $k$-d trees with rectangular cells, since the proximity of cells is handled by the MAC, and is unrelated to the tree structure. In other words, the dual tree traversal provides a simple bookkeeping strategy for constructing a MAC-based interaction stencil that is mutually exclusive at each level, by letting the target cells inherit a unique stack of source cells from their parents.

\section*{Results with the auto-tuning hybrid treecode/FMM method}
\vspace{\up}

\begin{figure*}
\centering
\includegraphics[width=\textwidth]{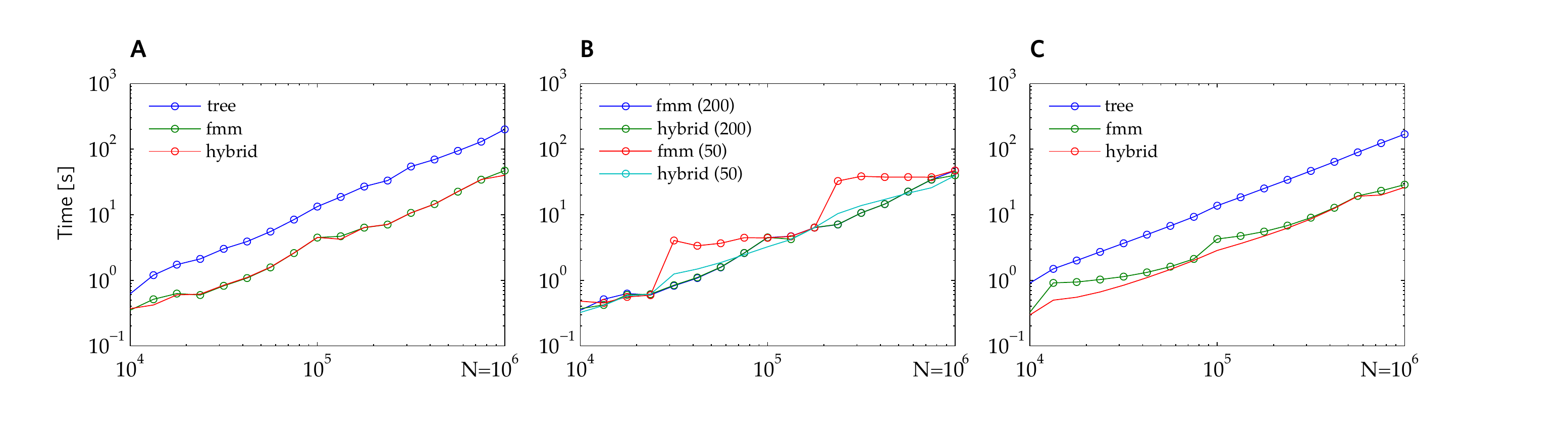}
\caption{Timings on CPU for Laplace kernel (potential+force), \textcolor{update}{treecode with $p=5$  and $p=8$ for FMM and hybrid}---(A) treecode, FMM and hybrid method, with particles randomly placed in a cube; (B) FMM and hybrid method with different values of $N_{crit}$; (C) particles randomly placed on a spherical shell.}
\label{fig:cpu}
\end{figure*}

The hardware used for these experiments was: Intel Xeon X5650 2.67GHz CPU, and NVIDIA \textcolor{update}{GeForce GTX590} GPU. In most cases, particles are randomly distributed in a cube of size $[-1,1]^3$ with the number of particles in the range \textcolor{update}{$N=10^4-10^6$ for the CPU runs, and $N=10^5-10^7$ for the GPU runs.} The calculations were performed for the Laplace kernel potential and force. Treecode, FMM, and hybrid method all use the same adaptive tree structure, dual tree traversals, and MAC-based interaction lists. \textcolor{update}{We define the MAC as $\theta=(r_t+r_s)/R$, where $r_t$ and $r_s$ are the radius of the target cell and source cell, respectively, and $R$ is the distance between the cell centers. We set $\theta=0.5$, which is equivalent to a typical FMM with a $3\times3\times3$ neighbor list.}
The only difference between treecode, FMM, and hybrid method is that \textcolor{update}{the  treecode always performs cell-particle interactions, and the FMM always performs cell-cell interactions,} while the hybrid can choose between cell-cell, cell-particle, and particle-particle interactions. All methods including the treecode use spherical harmonic expansions, and auto-tuning was not applied for the selection between different fast translation schemes at this time. Auto-tuning was only used to optimize the selection between, cell-cell, cell-particle, and particle-particle interactions in the hybrid method.

The timings on a single CPU core for treecode, FMM, and hybrid method are shown in Figure \ref{fig:cpu}(A). \textcolor{update}{The order of expansions is $p=5$ for the treecode and $p=8$ for the FMM, which yields an accuracy of 4 significant digits for the force. The results indicate that the hybrid method is always favoring cell-cell interactions and does not provide a visible advantage over the pure FMM. This seems to contradict previous works arguing that a mechanism to choose between cell-cell and cell-particle interactions should be optimal\cite{ChengETal1999}.} But since the treecode is implemented here using spherical expansions, the performance at low accuracy may be suboptimal compared to highly tuned Cartesian treecodes. \textcolor{update}{Thus, we conclude that a hybrid method that is faster than pure FMM may have to include more flexibility than just the dynamic choice of the type of interactions.}

\textcolor{update}{In treecodes and FMMs, the maximum number of particles per cell $N_{crit}$ is set (by the user) such that the loads of near-field evaluation and far-field evaluation are balanced. Setting this number to be too small will result in a very deep tree structure with a disproportionately large amount of far-field and too few near-field evaluations. Conversely, if $N_{crit}$ is set to be too large, the resulting tree structure will be too shallow and a large amount of time will be spent in the near-field evaluation. The advantage of the hybrid method becomes clear when we compare for two different values of $N_{crit}$ in Figure \ref{fig:cpu}(B). The first two legend entries are identical to those in Figure \ref{fig:cpu}(A), where a well-chosen value $N_{crit}=200$ was used. The latter two entries are the same tests, but with $N_{crit}=50$. In this case, the FMM suffers from load imbalance between the near-field and far-field evaluations, while the hybrid method does not because it can choose to perform particle-particle interactions even if the cell is not at the leaf-level.}

\textcolor{update}{To investigate the effect of adaptive distributions of particles, we performed additional tests with a random distribution of particles placed on a spherical shell. Figure \ref{fig:cpu}(C) shows the calculation time against the number of particles for this case. The maximum number of particles per cell was set to $N_{crit}=20$, a value that was as close to the optimum as we could get by manual adjustment, yet the hybrid method produced a better result by automatically fine-tuning the balance between the particle-particle and cell-cell interactions throughout the adaptive tree. We monitored the number of cell pairs which performed particle-particle, cell-particle, and cell-cell interactions. At $N=10^5$, the treecode executed 5 times more cell-particle interactions than particle-particle interactions, the FMM executed 3 times more cell-cell interactions than particle-particle interactions, and the hybrid executed two times more particle-particle interactions than cell-cell interactions. In conclusion, the cell-particle interaction never seems to give an advantage on CPUs, which means that the FMM is always faster than the treecode on CPUs for any given accuracy. The common belief that treecodes are faster for low accuracy must stem from the fact that typical implementations use Cartesian expansions for the treecode and spherical expansions for the FMM, and thus has nothing to do with the choice of cell-particle vs. cell-cell operations.}

\begin{figure}
\centering
\includegraphics[width=0.50\textwidth]{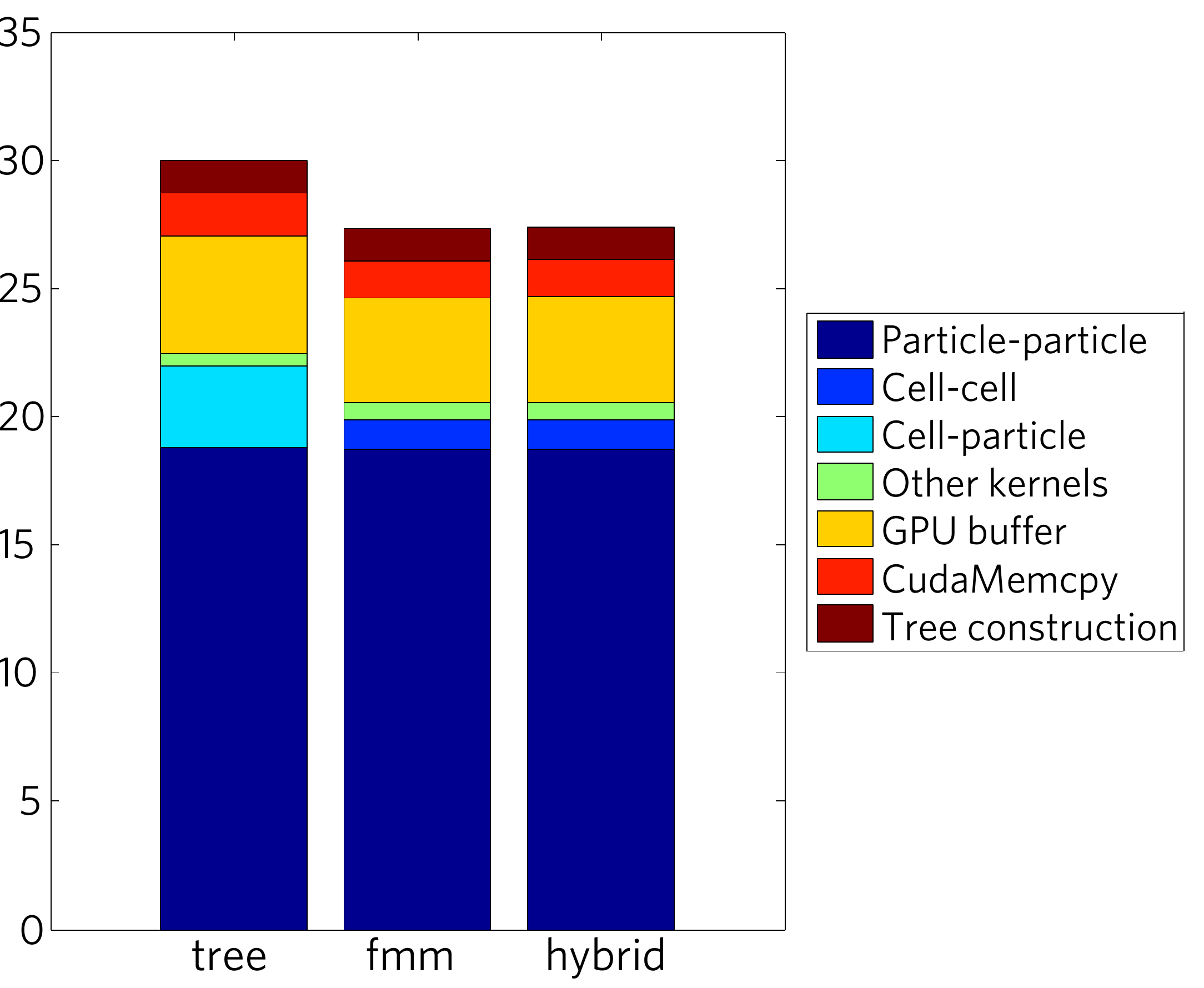}
\caption{Breakdown of the run time on GPU for the three methods. $N=10^7$ particles randomly placed in a cube, Laplace kernel potential+force, \textcolor{update}{$p=5$ for the treecode and $p=8$ for FMM and hybrid.} }
\label{fig:breakdown}
\end{figure}

\textcolor{update}{In similar experiments on GPU, all interactions are computed on the device (in single-precision) but tree construction is done on the CPU. Unlike the CPU case, where the treecode was significantly slower, the run times of the three method were very similar on GPU, and an equivalent plot to Figure \ref{fig:cpu}(C) is uninteresting.}
\textcolor{update}{Figure \ref{fig:breakdown} shows the breakdown of the calculation time on GPU for $N=10^7$, revealing that a large proportion of the time is being spent on particle-particle interactions. On GPU, it is relatively more costly to perform cell-cell interactions and shifting more work to the particle-particle interactions results in shorter runtime. Furthermore, auto-tuning on GPUs} faces the following problem: the calculation time is \emph{not} proportional to the problem size because larger problems are able to utilize more threads, and thus run more efficiently on GPUs. This makes it very difficult to predict the execution time of each kernel \textcolor{update}{by running a small test in the beginning.} Hence, the optimization between cell-cell, cell-particle, and particle-particle kernels may not function properly.

\begin{figure*}
\centering
\includegraphics[width=0.8\textwidth]{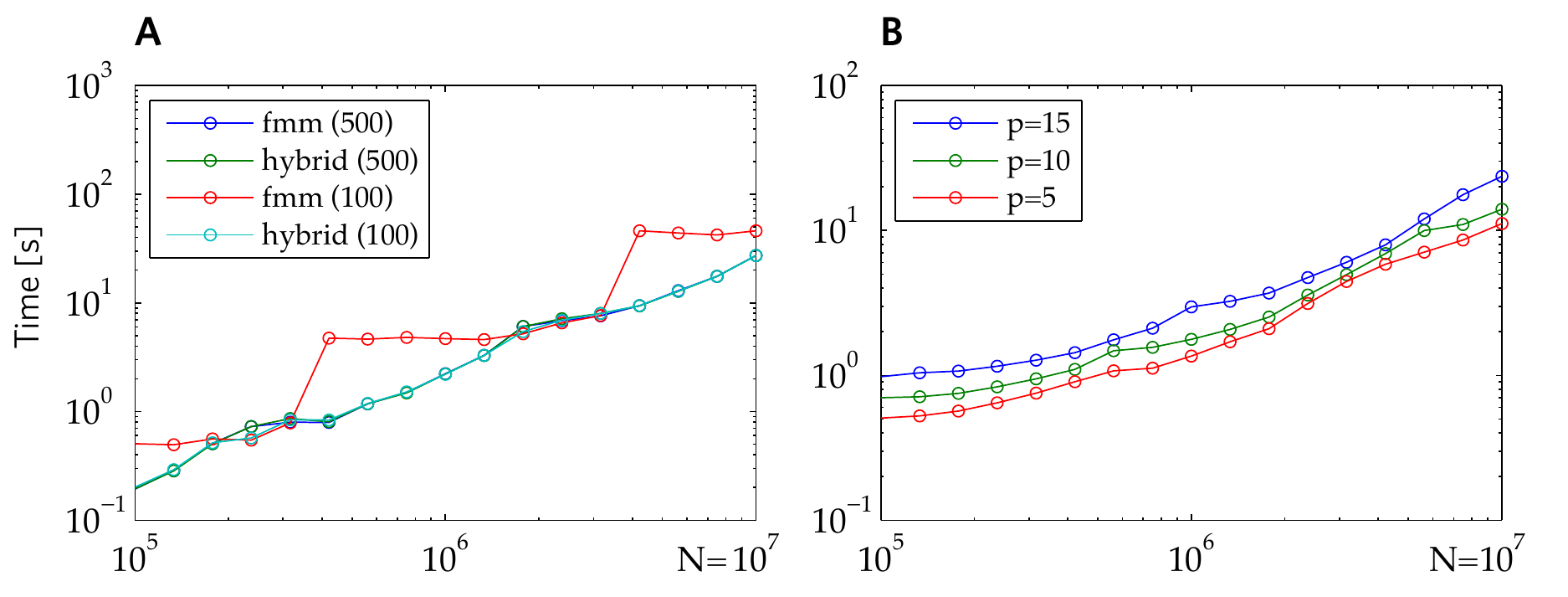}
\caption{Timings on GPU for Laplace kernel (potential+force) with particles randomly placed in a cube---(A) FMM and hybrid method for different values of $N_{crit}$, $p=8$; (B) hybrid method using different orders of expansion $p$.}
\label{fig:gpu2}
\end{figure*}

\textcolor{update}{In order to confirm that the auto-tuning capability is functioning on GPUs, we performed, as before, a test with different values of the maximum number of particles per cell, $N_{crit}$. As seen in Figure \ref{fig:gpu2}(A), for a choice of $N_{crit}=100$ the FMM experiences an imbalance between particle-particle and cell-cell interactions. In contrast, the hybrid method optimizes the balance between particle-particle and cell-cell interactions and achieves optimum performance for all $N$. We conclude that the auto-tuning mechanism is functioning on GPUs.}

The accuracy dependence of our hybrid treecode-FMM is shown in Figure \ref{fig:gpu2}(B). The calculation conditions are identical to the previous calculations except the order of expansion is changed from $p=5$ to $p=15$. The $p$-dependence is rather small considering the fact that we are using a $\mathcal{O}(p^4)$ cell-cell interaction kernel. One reason for this is that at the considered range of $p$ a smaller constant in front of the $p^4$ term allows the lower order terms to dominate. Another reason is the GPU being able to process the kernels with larger $p$ more efficiently, because they have a higher flop/byte rate. Our code is able to calculate the Laplace potential and force for $N=10^7$ particles with $p=15$ in approximately $23$ seconds on a single GPU.

\bigskip

\lettrine{\textcolor{dropped}{W}}{} ith the current hybridization of treecode and FMM, combined with auto-tuning capabilities on heterogeneous architectures, the flexibility of fast $N$-body methods has been greatly enhanced. The fact that the current method can automatically choose the optimal interactions, on a given heterogeneous system, alleviates the user from two major burdens. \textcolor{update}{Firstly, the user does not need to decide among treecode or FMM, predicting which algorithm will be faster for a particular application given the accuracy requirements---they are now one algorithm.}  Secondly, there is no need to tweak parameters, e.g., particles per cell, in order to achieve optimal performance on GPUs---the same code can run on any machine without changing anything. This feature is a requirement to developing a black-box software library for fast $N$-body algorithms \emph{on heterogeneous systems}, which is our goal. 

Our codes are available for unrestricted use, under the MIT license; to obtain the codes and run the tests in this paper, the reader may follow instructions in the website at \href{http://www.bu.edu/exafmm/}{\textbf{www.bu.edu/exafmm}}.

\section*{Acknowledgements}
\vspace{\up}

{\sf \emph{We're grateful for the support from the US National Science Foundation and the Office of Naval Research. Recent grant numbers are NSF OCI-0946441, and ONR award \#N00014-11-1-0356.}

\bigskip

}

\bibliographystyle{unsrt}
\bibliography{FastMethods,scicomp}

\vspace{1cm}

\small
{\sf

\noindent \textbf{Rio Yokota} obtained his PhD from Keio University in 2009 and went on to work as a postdoctoral researcher with Prof.\ Barba at the University of Bristol and Boston University.   During his PhD, he worked on the implementation of fast $N$-body algorithms on special purpose machines such as \textsc{mdgrape}-3, and then on GPUs after CUDA was released.  He co-authored a paper awarded the Gordon Bell prize in 2009 in the price/performance category, using GPUs. Since Fall 2011, he joins King Abdullah University of Science and Technology as a research scientist in the Center for Extreme Computing.

\bigskip

\noindent\textbf{Lorena A. Barba} is an Assistant Professor of Mechanical Engineering at Boston University.  She obtained her PhD in Aeronautics from the California Institute of Technology in 2004, and then joined Bristol University, UK, as a Lecturer in Applied Mathematics.  Her research interests include computational fluid dynamics, especially particle methods for fluid simulation and immersed boundary methods; fundamental and applied aspects of fluid dynamics, especially flows dominated by vorticity dynamics; the fast multipole method and applications; and scientific computing on GPU architecture. She received the Amelia Earhart Fellowship, a First Grant award from the UK Engineering and Physical Sciences (EPSRC), and an NVIDIA Academic Partnership grant in August 2011.
}

\end{document}